\documentclass[11pt]{amsart}
\usepackage{amsmath,amsfonts,amsthm,amscd,amssymb,mathrsfs,amssymb,bm}
\usepackage{mathrsfs}
\usepackage{graphicx}
\usepackage[all]{xy}

\newtheorem{theorem}{Theorem}

\newtheorem{lemma}[theorem]{Lemma}

\newtheorem{remark}[theorem]{Remark}
\newcommand{\CP}{\mathbb{CP}}

\newcommand{\ol}{\overline}
\newcommand{\lra}{\longrightarrow}

\newcommand{\set}{\,|\,}

\numberwithin{equation}{section}
\numberwithin{theorem}{section}
\begin{document}
\bibliographystyle{alpha} 
\title{Degenerations of LeBrun twistor spaces}
\author{Nobuhiro Honda}
\address{Department of Mathematics, Tokyo Institute of Technology, O-okayama,  Tokyo, Japan}
\address{Current Address: Mathematical Institute, Tohoku University, Sendai, Miyagi, Japan}
\email{honda@math.tohoku.ac.jp}
\thanks{The author has been partially supported by the Grant-in-Aid for Young Scientists  (B), The Ministry of Education, Culture, Sports, Science and Technology, Japan.}
\begin{abstract}
We investigate various limits of the twistor spaces associated to the self-dual metrics on $n \mathbb{CP}^2$, the connected sum of the complex projective planes, constructed by C. LeBrun.
In particular, we explicitly present the following 3 kinds of degenerations
whose limits of the corresponding metrics are:
(a) LeBrun metrics on $(n-1) \mathbb{CP}^2$,
(b) (another)  LeBrun metrics on the total space of the line bundle $\mathscr O(-n)$ over $\mathbb{CP}^1$, 
(c) the hyper-K\"ahler metrics on the small resolution of rational double points of type $A_{n-1}$, constructed by G.W. Gibbons and S.W. Hawking.
\end{abstract}
\maketitle
\setcounter{tocdepth}{1}
\vspace{-5mm}
\section{ Introduction.}

According to a decomposition of the space of the curvature tensors by  natural action of an orthogonal group, 
the curvature tensor of a Riemannian metric splits into a traceless Ricci tensor, scalar curvature, and the Weyl curvature tensor.
As is well-known, the Weyl tensor is invariant under conformal changes of the metric, and when the manifold is more than 3-dimensional,
the metric is conformally flat if and only if the Weyl tensor identically vanishes.
When the manifold is 4-dimensional and oriented, since the Lie algebra $\mathfrak{so}(4)$ splits as $\mathfrak{su}(2)\oplus\mathfrak{su}(2)$, the bundle of 2-forms, which is exactly the bundle associated to the adjoint representation of $\rm{SO}(4)$, splits into two components, the self-dual  and anti-self-dual parts.
Accordingly, the Weyl curvature tensor splits into  two parts.
A Riemannian metric on an oriented 4-manifold  is said to be {\em self-dual}\, if the anti-self-dual part of the Weyl tensor vanishes everywhere.

Self-dual metrics seem to have been paid special attentions  because 
of its deep connection with complex geometry. 
Namely, based on an idea by R. Penrose, it is shown by Atiyah, Hitchin and Singer \cite{ahs}  that the total space of a projectified spin bundle over any oriented Riemannian 4-manifold carries a natural almost complex structure, whose integrability condition is precisely the self-duality of the metric.
Thus any self-dual metric canonically produces a 3-dimensional complex manifold; this is what is called the {\em twistor space\,}
associated to a self-dual metric.

While a considerable number of self-dual metrics and twistor spaces are  known now, 
 the most  illuminating and tractable one  still seems to be the LeBrun metrics and their twistor spaces \cite{LB91}.
The LeBrun metrics are parameterized by a finite number of points (called {\em monopole points}) on the hyperbolic 3-space; more precisely, the self-dual metrics are constructed on the total space of certain U(1)-bundle over the hyperbolic 3-space with  $n$  monopole points removed, and
this metric conformally extends to a compactification which is exactly $n\mathbb{CP}^2$.  
The two self-dual manifolds determined from two configurations of $n$ monopole points are conformally equivalent if and only if the configurations are equivalent under the usual PSL$(2,\mathbb C)$-action on the hyperbolic space.

By moving the monopole points on the hyperbolic space into extreme configurations, it is possible to consider various natural interesting limits of the LeBrun metrics.
When such a deformation (to an extreme configuration) is explicitly given, often it looks not so difficult to  identify the limit self-dual structure, at least heuristically.
However, it seems not obvious at all (at least for the author) what happens for the twistor spaces at the limit.
A purpose of this paper is to give some explicit and rigorous answers to this question.  

In Section 2, we consider the case where one of the $n$  monopole points on the hyperbolic space goes to a point at the ideal boundary.
In this case, the limit space is  $(n-1)\mathbb{CP}^2$ equipped with a LeBrun metric. In particular, the topology of the 4-manifold changes.
At the level of twistor spaces, we shall show the following.
First if we take the corresponding natural limit for the twistor spaces on $n\mathbb{CP}^2$, then one obtains a 3-fold having one ordinary double point.
In this degeneration, certain degree-one divisor (and its conjugation) breaks into two irreducible components, exactly one of which is $\mathbb{CP}^2$.
(So together with its conjugation, we have two $\mathbb{CP}^2$-s.
The intersection of these is the double point.)
If we take an appropriate small resolution for this ordinary double point, then the normal bundle of the two $\mathbb{CP}^2$-s becomes isomorphic to $\mathscr O(-1)$.
If we blow-down these, then we obtain a LeBrun  twistor space on $(n-1)\mathbb{CP}^2$.
 This is the main result in Section 2.
 We remark that our actual explanation is a bit different in that we mainly use {\em projective models} of the twistor spaces and not the twistor spaces themselves.
 (Since the necessary operation for obtaining the twistor spaces from projective models are given by LeBrun, this is not a serious matter.)
The use of projective models is indispensable when taking a limit.

In Section 3, differently from the limit in Section 2, we investigate the case where all the $n$ points on the hyperbolic space approach to a same point.
In this case we show that the limit of the twistor space is the twistor space of LeBrun's scalar-flat K\"ahler metric on the total space of the line bundle $\mathscr O(-n)\to \mathbb{CP}^1$.
The main body of the proof of this result is to find out all twistor lines in the limit space, to verify that the space is actually a twistor space, and we finally use a Pontecorvo's result to a certain divisor and a LeBrun's characterization result of the last metric, in order to identify the self-dual structure for the twistor space.

In Section 4, we consider another kind of limits.
Namely, instead of moving the monopole points, we fix these points and take rescaling to the hyperbolic space which makes the curvature zero.
Then in the limit the hyperbolic space becomes the flat Euclidean space and the LeBrun metric approaches to the Gibbons-Hawking's hyperK\"ahler metric on a 4-manifold diffeomorphic to $\mathbb C^2/\Gamma$, where $\Gamma$ is a cyclic subgroup of SU$(2)$ of order $n$.
Note that the twistor spaces for the last metrics are explicitly constructed by Hitchin \cite{Hi79}.
For any LeBrun twistor space on $n\mathbb{CP}^2$, we construct an explicit  degenerating 1-dimensional family such that general fibers are  Zariski-open subsets of LeBrun twistor spaces, and such that the limit fiber is the above twistor space of Hitchin's.
This completely agrees with the situation for the metrics
which was obtained in \cite[Section 5]{LB91}.
Finally we see that if we pull-back the above degenerating family by certain trivializing map outside the origin, then over the central fiber there appears the twistor space of the LeBrun's scalar-flat K\"ahler metric on $\mathscr O(-n)$.

We should mention what is happening for the base 4-manifolds for these degenerations.
In the final section we give a brief discussion about convergence of the LeBrun metrics, focusing on their  K\"ahler representatives. 
For investigation on various limits of LeBrun metrics
 from differential geometric point of view, 
we refer a recent paper by Jeff Viaclovsky \cite{V}.

\vspace{3mm} 
\noindent{\bf Acknowledgement.}
First of all, I would like to express my sincere gratitude to Jeff Viaclovsky for a number of helpful discussions about the material in this paper. 
Certainly, without these discussions, I could not get any result in this paper.
Also I would like to thank Claude LeBrun for  suggesting the use of minitwistor spaces for connecting his twistor spaces with Hitchin's one.
This suggestion was crucial for obtaining the result in Section 4. 
I also would like to thank the referees for a lot of kind and insightful suggestions, which yields Section 5.
At this point I also would like to thank Kazuo Akutagawa for teaching materials on convergence of critical metrics on 4-manifolds, and kindly giving me a lot of suggestions.

\section{Degenerations to LeBrun metrics on $(n-1)\mathbb{CP}^2$ from $n\mathbb{CP}^2$}\label{s:2}
\subsection{LeBrun twistor spaces}
In this section we briefly recall LeBrun's  construction of his twistor spaces on $n \mathbb{CP}^2$ \cite{LB91} . 
These are needed for giving our explicit construction of degenerations.

Put $Q:=\mathbb{CP}^1\times\mathbb{CP}^1$, equipped with the real structure $
\sigma:(u,v)\mapsto (\ol{v},\ol{u}),
$
where $u$ and $v$ are non-homogeneous coordinates on the two factors  respectively. 
For each  non-negative integer $n$ consider a rank-3 vector bundle over $Q$ defined by
\begin{equation}
\mathscr E_n:=\mathscr O(n-1,1)\oplus\mathscr O(1,n-1)\oplus \mathscr O.
\end{equation}
Then as $\sigma$ interchanges the two factors, there is a natural anti-linear isomorphism  between the line bundles  ${\mathscr O(n-1,1)}$ and $\mathscr O(1,n-1)$.
This induces a real structure on the $\mathbb{CP}^2$-bundle
$\mathbb P(\mathscr E_n)=(\mathscr E_n-\{0\})/\mathbb C^*$, for which we still denote by $\sigma$.


Next we take any different $n$ points on the upper-half space $\mathscr H^3$ and consider the LeBrun metrics on $n \mathbb{CP}^2$ determined from them.
The twistor space of this metric is constructed as follows.
Let $\mathscr C_1,\cdots,\mathscr C_n$ be the $(1,1)$-curves on $Q$ corresponding to the $n$ points.
Each of these curves are real and has no real point.
So any two different $\mathscr C_i$ and $\mathscr C_j$ intersect transversally at two points.
On the affine plane $\mathbb C^2\subset Q$ on which the coordinates $(u,v)$ are valid,
the curve $\mathscr C_i$ is defined by the equation
\begin{equation}\label{P_i}
P_i(u,v)=a_iuv+b_iu+c_iv+d_i,
\end{equation}
where by reality the coefficients satisfy $a_i,d_i\in\mathbb R$ and  $b_i=\ol{c}_i$.
Then define  an algebraic subvariety of $\mathbb P(\mathscr E_n)$ by
\begin{equation}\label{LB-n}
X:=\{(x,y,z)\in\mathbb P(\mathscr E_n)\set xy=P_1(u,v)P_2(u,v)\cdots P_n(u,v)\,z^2\},
\end{equation}
where $(x,y)\in\mathscr O(n-1,1)\oplus\mathscr O(1,n-1)$ and $z\in\mathscr O$.
We call $X$ {\em the projective model of a LeBrun twistor space on\, $n \mathbb{CP}^2$}.
 $X$ has an obvious structure of a conic bundle over $Q$.
As $xy\in\mathscr O(n,n)$ and $P_1\cdots P_n\in H^0(\mathscr O(n,n))$, 
there exists no discriminant curve other than $\mathscr C_i$-s. 
The singular locus of $X$ is precisely over the intersection points of $\mathscr C_i$ and $\mathscr C_j$ $(i\neq j)$ and they are lying on the section $\{x=y=0\}$.
From the equation in \eqref{LB-n}, the conic bundle 
 $X$ has distinguished, mutually conjugate two sections 
\begin{equation}\label{2sect}
E:=\{x=z=0\}\hspace{2mm}{\text{and}}\hspace{2mm}\ol{E}:=\{y=z=0\}.
\end{equation}
These are clearly away from the singular locus of $X$.
Their normal bundles are
\begin{equation}\label{nb}
N_{E/X}\simeq\mathscr O(-1,1-n), \hspace{2mm}
N_{\ol{E}/X}\simeq\mathscr O(1-n,-1),
\end{equation}
and therefore $E$ and $\ol{E}$ can be blown-down along one of the two projections $\mathbb{CP}^1\times\mathbb{CP}^1\to\mathbb{CP}^1$.
Let $X\to Y$ be these blowing-down of $E$ and $\ol{E}$.
Then the required LeBrun twistor space is obtained from $Y$ by taking appropriate small resolutions for the singularities over $\mathscr C_i\cap\mathscr C_j$.

Since any LeBrun metric admits a  U(1)-action, 
LeBrun twistor spaces  always admit a holomorphic $\mathbb C^*$-action.
It naturally induces a $\mathbb C^*$-action on  the projective model $X$ and in the above coordinates it is explicitly given by
\begin{equation}\label{action1}
(x,y,z)\longmapsto (sx,s^{-1}y,z),\,\,s\in\mathbb C^*.
\end{equation}

\subsection{A degeneration from $2\mathbb{CP}^2$ to $\mathbb{CP}^2$}
In this subsection we consider the case $n=2$.
In this case by normalizing the coordinates $(u,v)$, we may suppose that the projective model $X$ is defined by
\begin{equation}\label{LB-2}
xy=(uv+1)(uv+\lambda)\,z^2
\end{equation}
where $0<\lambda<1$.
The parameter $\lambda$ represents the conformal class of the metric.

From now on we investigate what happens when $\lambda\to 0$.
Namely we examine the variety
\begin{equation}\label{LBlim1}
X_{\infty}:=\{(x,y,z)\in\mathbb P(\mathscr E_2)\set xy=uv(1+uv)z^2\}.
\end{equation}
(The subscript $\infty$ reflects the fact that   this limit corresponds to the situation that one of the two  monopole points  goes to infinity.)
We name the two discriminant curves as
$$
\mathscr C_1^{\infty}:=\{uv=0\},\quad
\mathscr C_2:=\{uv+1=0\}
$$
Of course, the reducible curve $\mathscr C_1^{\infty}$ is a limit of the irreducible discriminant curve $\{uv+\lambda=0\}$. 
The singularities of the curve $\mathscr C_1^{\infty}+\mathscr C_2$ consists of 3 points $(u,v)=(0,0), (0,\infty)$ and $(\infty,0)$.
Correspondingly  $X_{\infty}$ has 3 ordinary double points over there.
The divisors $E$ and $\ol{E}$ are still contained in $X_{\infty}$ and disjoint from these singularities,  satisfying \eqref{nb}.

We are going to show that the variety $X_{\infty}$ is explicitly birational to the flag twistor space of the Fubini-Study metric on $\mathbb{CP}^2$ along the following two steps:

(i) Choose appropriate small resolutions of all singularities of $X_{\infty}$,

(ii) Blow-down  two redundant divisors into curves.

\begin{figure}
\includegraphics{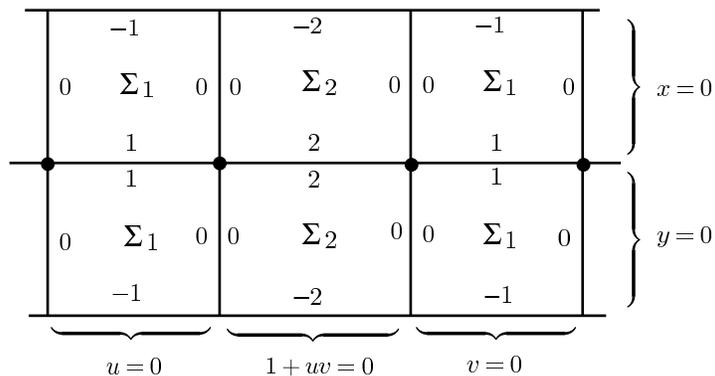}
\caption{The inverse image of $\mathscr C_1^{\infty}+\mathscr C_2$ under $X_{\infty}\to Q$.
(The left and the right vertical edges represents the same fiber over $(0,0)$.)}
\label{fig-2_1}
\end{figure}

\begin{figure}
\includegraphics{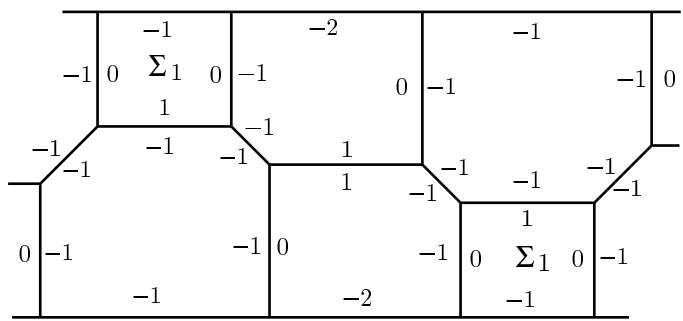}
\caption{After taking small resolutions. (The left three edges and the right ones are identified as in the case for Figure \ref{fig-2_1})}
\label{fig-2_2}
\end{figure}

\noindent
Figure \ref{fig-2_1} displays  the inverse image of the discriminant curve $\mathscr C_1^{\infty}+\mathscr C_2$.
Each face (square) represents an irreducible component of the inverse image of  $\mathscr C_1^{\infty}$
or $\mathscr C_2$, and
the numbers attached to each edges represent the intersection number of the edge (= a rational curve) in the component.
The dotted points (which are exactly the points where 4 faces meet) are the ordinary double points of $X_{\infty}$.
Further $\Sigma_k$ means the Hirzebruch surface $\mathbb P(\mathscr O(k)\oplus\mathscr O)$.
Then as the first step (i) we take  small resolutions of the ordinary double points which transform the divisors in Figure \ref{fig-2_1} into the situation displayed in Figure \ref{fig-2_2}.
The resulting 3-fold is non-singular.
In this 3-fold, there remains exactly two divisors which are isomorphic to $\Sigma_1$.
Since their normal bundles are degree $(-1)$ along fibers as  in Figure \ref{fig-2_2}, 
 they can be blown down in that direction.
Let $\mathbb F'$ be the resulting threefold.
Again $\mathbb F'$ is non-singular and has a conic bundle structure $\mathbb F'\to Q$ induced from that on $X_{\infty}$.
The discriminant locus of this is a single irreducible curve $\mathscr C_2$.
Namely, both of the two irreducible components of $\mathscr C_1^{\infty}$ are {\em not} discriminant curve any more.
Thus  our modification (especially the blow-down step (ii)) has an effect of {\em removing the reducible curve $\mathscr C_1^{\infty}$ from the discriminant locus.}

\begin{figure}
\includegraphics{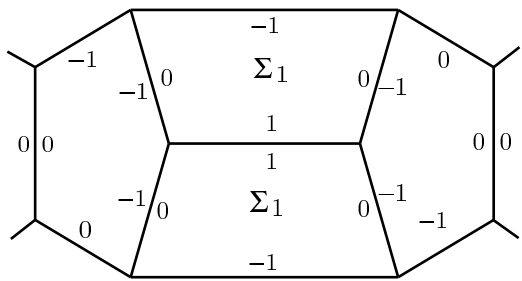}
\caption{The inverse image of $\mathscr C_1^{\infty}+\mathscr C_2$ under $\mathbb F'\to Q$.
(The left and right vertical edges are identified.)}
\label{fig-2_3}
\end{figure}

On the other hand, since  LeBrun metric on $\mathbb{CP}^2$ is nothing but the Fubini-Study metric,  the flag twistor space $\mathbb F$ can be obtained from the `projective model'
$$
\tilde{\mathbb F}:=\{(x,y,z)\in\mathbb P(\mathscr E_1)\set xy=(1+uv)z^2\}
$$
 by blowing down two sections $E$ and $\ol{E}$.
Obviously the two conic bundles $\mathbb F'\to Q$ and $\tilde{\mathbb F}\to Q$  have the same discriminant locus $\mathscr C_1$.
 More strongly, we have
 
\begin{theorem}\label{thm-2}
$\mathbb F'$ and  \,$\tilde{\mathbb F}$ are biholomorphic.
\end{theorem}

This implies that one can obtain the flag twistor space on $\mathbb{CP}^2$ from a LeBrun twistor space on $2 \mathbb{CP}^2$ by first taking a limit ($\lambda\to 0$) and then blowing down redundant divisors.

We postpone a proof of  Theorem \ref{thm-2}  until next subsection in which we prove a theorem which contains 
 Theorem \ref{thm-2} as a special case.

\subsection{A degeneration from $n\mathbb{CP}^2$ to $(n-1) \mathbb{CP}^2$}
In this section we generalize the construction in the previous subsection to the case of $n \mathbb{CP}^2$  and prove that  LeBrun twistor spaces on $(n-1) \mathbb{CP}^2$ can be obtained from those on $n \mathbb{CP}^2$ by first taking a  limit and then applying some  birational transformations.

\begin{figure}
\includegraphics{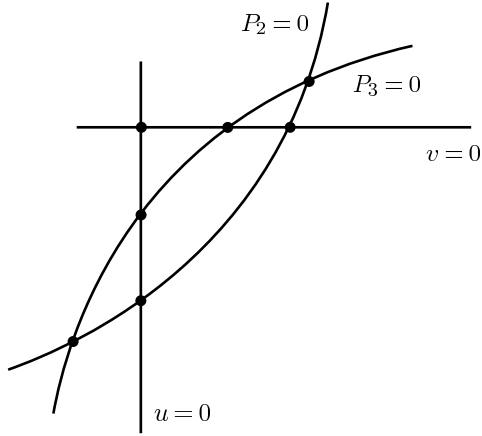}
\caption{The discriminant locus of  $X_{\infty}\to Q$ for the case $n=3$.}
\label{fig-n_0}
\end{figure}

Suppose $n\ge 3$. 
We keep notation in the last subsection.
First we take appropriate  coordinates on $Q=\mathbb{CP}^1\times\mathbb{CP}^1$ to normalize defining equation \eqref{LB-n}  as
\begin{equation}
xy=(uv+\lambda)P_2(u,v)P_3(u,v)\cdots P_n(u,v)z^2,\quad\lambda\in\mathbb R^{\times}.
\end{equation}
(We can further suppose $P_2(u,v)=uv+1$, but we do not assume this because the curve $P_2=0$ does not play a particular role.)
Without loss of generality we can suppose that if $i\neq j$ the curves $\{P_i=0\}$ and $\{P_j=0\}$ do not intersect on $\{uv=0\}$.
If we put $\lambda=0$, we obtain a variety
\begin{equation}\label{limLB-n-2}
X_{\infty}:=\{(x,y,z)\in \mathbb P(\mathscr E_n)\set xy=uvP_2(u,v)P_3(u,v)\cdots P_n(u,v)z^2\}.
\end{equation}
Then again all singularities of $X_{\infty}$ are over singular points of the discriminant curves 
$$uvP_2(u,v)P_3(u,v)\cdots P_n(u,v)=0,$$
 and they are also lying on the section $\{x=y=0\}$.
Generically all of them are ordinary double points.
(But we do not suppose it.)
In Figure \ref{fig-n_0} the discriminant locus is illustrated, in the case $n=3$.

\begin{figure}
\includegraphics{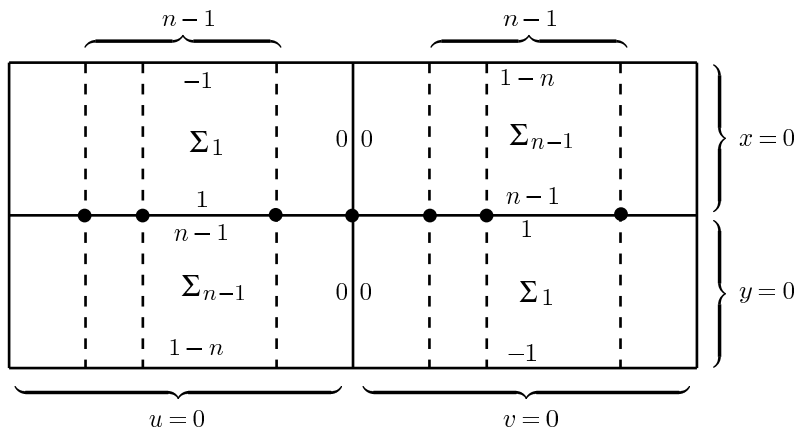}
\caption{The inverse image of the reducible curve $\{uv=0\}$ under $X_{\infty}\to Q$}
\label{fig-n_1}
\end{figure}

As in the case $n=2$, we consider the inverse image of the reducible curve $\{uv=0\}$. 
The inverse image of an irreducible component $\{u=0\}$ consists of two non-singular divisors, which are biholomorphic to $\Sigma_1$ and $\Sigma_{n-1}$ respectively.
The situation is similar to the curve $\{v=0\}$, although the place of  $\Sigma_1$ and $\Sigma_{n-1}$ are upside down as illustrated in Figure \ref{fig-n_1}.
Thus there are exactly two $\Sigma_1$-s in the inverse image of the curve $\{uv=0\}$.
On each of these two $\Sigma_1$-s there are precisely $n$ ordinary double points of $X_{\infty}$.
Exactly one of these points are shared by the two $\Sigma_1$-s.
This point is over the point $\{u=v=0\}$.
In Figure \ref{fig-n_1} these $(2n-1)$ singularities are displayed as dotted points.
(The last ordinary point is at the center.)
The vertical broken lines represent the intersection with the inverse image of the discriminant curve $\{P_i=0\}, \,2\le i\le n$.

\begin{figure}
\includegraphics{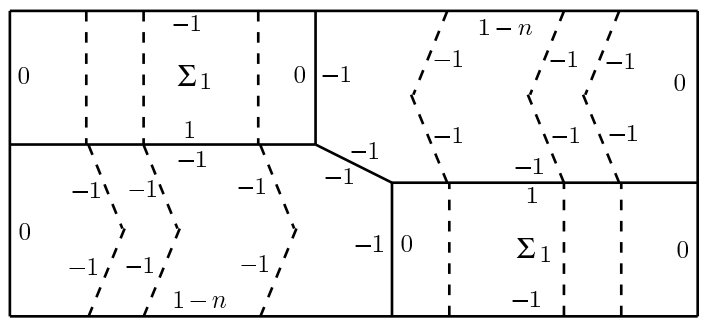}
\caption{The inverse image of $\{uv=0\}$ under $Y\to Q$
(i.e.\,\,strict transforms of the divisors in Figure \ref{fig-n_1}).}
\label{fig-n_2}
\end{figure}

Then as the first step of the present modification, we resolve all these $(2n-1)$ ordinary double points by small resolutions which are uniquely determined  by the rule that {\em the two $\Sigma_1$-s are always unchanged}.
Similarly to the case $n=2$,
we denote the resulting space by $Y$.
There are still other  double points of $Y$  over the intersections of $\{P_i=0\}$ and $\{P_j=0\}$ for $2\le i<j\le n$.
(These double points did not appear in the case $n=2$.)
We do not touch these singularities.
In $Y$, the inverse image of the curve $\{uv=0\}$ is as illustrated in Figure \ref{fig-n_2}.
In particular $Y$ is non-singular in a neighborhood of these divisors.
Further as in the figure  the normal bundles of the two $\Sigma_1$-s in $Y$ is degree $(-1)$ along fibers of $\Sigma_1\to\mathbb{CP}^1$.
Therefore these two divisors can be blown-down.
Let $Y\to X'$ be this blowing down.
This is the second step of the present modification.
In this way starting from a LeBrun twistor space on $n \mathbb{CP}^2$ we obtained a new space $X'$.
Let $f:X'\to Q$ be the projection naturally induced from the conic bundle projection $X_{\infty}\to Q$.
Then as in the case $n=2$, the curve $\{uv=0\}$ is {\em not} a discriminant curve of $f$.
\begin{figure}
\includegraphics{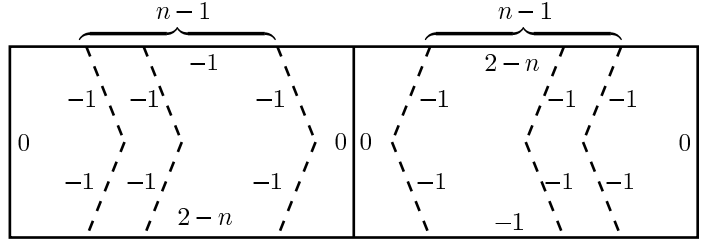}
\caption{The inverse image of $\{uv=0\}$ under $X'\to Q$ (and $\tilde X\to Q$).
(i.e.\,\,Transformations of the divisors in Figure \ref{fig-n_2}.)}
\label{fig-n_3}
\end{figure}
Then we have the following

\begin{theorem}\label{thm-3}
The above constructed threefold \,$X'$ is biholomorphic to the projective model (see (2.1)) of the LeBrun twistor space on $(n-1) \mathbb{CP}^2$ defined by
\begin{equation}\label{LB-(n-1)}
\{(x,y,z)\in \mathbb P(\mathscr E_{n-1})\set xy=P_2(u,v)\cdots P_n(u,v)z^2\}.
\end{equation}
\end{theorem}

\begin{proof}
Recall that the original conic bundle $X_{\infty}\to Q$ has two distinguished sections $E$ and $\ol E$, and 
the birational transformations from $X_{\infty}$ to $X'$ do not contract these sections.
So the projection $f:X'\to Q$ still has two distinguished sections, still denoted by $E$ and $\ol{E}$.
As explained in the preceding paragraph of Theorem \ref{thm-3}, $X'$ has  double points which are over $P_i=P_j=0$ for $2\le i<j\le n$.
These are lying on the section  $\{x=y=0\}$ and hence  disjoint from $E$ and $\ol{E}$.
We take any small resolutions of all these  double points of $X'$ to obtain a non-singular threefold $\tilde{X}$. 
Then as the curve $\{uv=0\}$ is not a discriminant curve for $f$, the discriminant locus of the composition $\tilde f:\tilde X\to X'\to Q$ is the curve $\{P_2(u,v)\cdots P_n(u,v)=0\}$.
Any fiber over  the discriminant locus is a chain of  smooth rational  curves.
Since the small resolution $\tilde X\to X'$ does not touch, $\tilde f$ still has distinguished sections  $E$ and $\ol{E}$.
As one can read off from  Figure \ref{fig-n_3}
(the upper horizontal line is the intersection with $E$ (or $\ol E$), and the lower one is the intersection with $\ol E$ (or $ E$)), their normal bundles in $\tilde X$ (or in $X'$)  are given by
\begin{equation}\label{nb-3}
N_{E/\tilde X}\simeq \mathscr O(-1,2-n),\,\,\,\text{and}\,\,\,N_{\ol{E}/\tilde X}\simeq\mathscr O(2-n,-1).
\end{equation}
Hence considering an exact sequence
\begin{equation}\label{ses1}
0\lra
\mathscr O_{\tilde X}\lra
\mathscr O_{\tilde X}(E+\ol{E})\lra
N_{E/\tilde X}\oplus N_{\ol{E}/\tilde X}\lra 0,
\end{equation}
 taking a direct image by $\tilde f$, noting $R^1\tilde f_*\mathscr O_{\tilde X}=0$ (as all fibers of $\tilde f$ is a chain of smooth rational curves), and computing Ext$^1$ for this sequence, we obtain an isomorphism
\begin{equation}
\tilde f_*\mathscr O_{\tilde X}(E+\ol{E})
\simeq
\mathscr O\oplus \mathscr O(-1,2-n)\oplus\mathscr O(2-n,-1).
\end{equation}
In particular, the direct image $f_*(\mathscr O_{\tilde X}(E+\ol{E}))$ is isomorphic to the dual bundle  $\mathscr E^*_{n-1}$, in which the projective model of a LeBrun twistor space on $(n-1) \mathbb{CP}^2$ is embedded.
By the morphism $\tilde f:\tilde X\to Q$ and the line bundle $\mathscr O_{\tilde X}(E+\ol{E})$ we obtain a rational map $\Phi:\tilde X\to\mathbb P(\mathscr E_{n-1}^*)^*=\mathbb P(\mathscr E_{n-1})$ over $Q$.
The restriction of $\Phi$ to a fiber $\tilde{f}^{-1}(p)$ is precisely the rational map associated to the 
linear system of the restricted line bundle $\mathscr O_{\tilde X}(E+\ol{E})|{\tilde{f}^{-1}(p)}$. 
It follows that $\Phi|\tilde f^{-1}(p)$ is an embedding  if $p\not\in \{P_i=P_j=0\}, \,2\le i<j\le n$. 
Further $\Phi$ contracts the exceptional curves of the small resolution $\tilde X\to X'$ since  $E$ and $\ol{E}$ are disjoint from the exceptional curves.
This means that $\Phi$ induces an isomorphism from $X'$ to $\Phi(\tilde X)$.
Namely the variety $X'$ can be embedded in $\mathbb P(\mathscr E_{n-1})$ as a conic sub-bundle.
It is obvious that the discriminant locus of the projection $X'=\Phi(\tilde X)\to Q$ is still the curve  $\{P_2(u,v)\cdots P_n(u,v)=0\}$
and the inverse image of its irreducible component $\{P_i=0\}$ $(2\le i\le n$) consists of two irreducible components.
This is the same for the projective model \eqref{LB-(n-1)}.
Using this fact and the $\mathbb C^*$-action \eqref{action1} we now show that $X'=\Phi(\tilde X)$ is biholomorphic  to the threefold \eqref{LB-(n-1)}.

Since the two sections $E$ and $\ol{E}$ in $X_{\infty}$ are fixed by  the $\mathbb C^*$-action \eqref{action1}, and the birational transformations from $X_{\infty}$ to $\tilde X$ clearly preserve this action, $E$ and $\ol{E}$ are still fixed under the natural $\mathbb C^*$-action on $\tilde X$.
Hence the sequence \eqref{ses1} is $\mathbb C^*$-equivariant and the morphism $\Phi:\tilde X\to\mathbb P(\mathscr E_{n-1})$ is $\mathbb C^*$-equivariant.
If we note that the $\mathbb C^*$-action on $N_{E/X_{\infty}}$ and $N_{\ol{E}/X_{\infty}}$ are given by scalar multiplications of $s$ and $s^{-1}$ respectively for $s\in\mathbb C^*$, 
it follows that naturally induced $\mathbb C^*$-action on $\mathscr E_{n-1}=\tilde f_*\mathscr O_{\tilde X}(E+\ol{E})$ is given by
$(x,y,z)\mapsto (sx,s^{-1}y,z)$.
Hence since the image $\Phi(\tilde X)=X'\subset\mathbb P(\mathscr E_{n-1})$ is $\mathbb C^*$-invariant, the equation of $\Phi(\tilde X)$ must be of the form
\begin{equation}\label{quad4}
xy=g(u,v)z^2
\end{equation}
for some $g(u,v)\in H^0(\mathscr O(n-1,n-1))$.
(The coefficients of the LHS must be constant by equi-dimensionality of $X'\to Q$.)
Then since the discriminant locus of  $X'\to Q$ is $\{P_2\cdots P_n=0\}$  we obtain $g(u,v)=cP_2(u,v)\cdots P_n(u,v)$ for some $c\in\mathbb C^*$.
Finally, by reality, we obtain $c\in\mathbb R^{\times}$ and hence
by normalizing the coordinate $z$, we can suppose $c=1$.
Thus we have seen that a defining equation of $X'$ in $\mathbb P(\mathscr E_{n-1})$ is the same as \eqref{LB-(n-1)}.
Hence $X'$ is isomorphic to \eqref{LB-(n-1)}.
\end{proof}

\subsection{Remarks about the construction}
The series of operations above starts from  a {\em projective model} of LeBrun twistor space, and  also ends at a projective model.
In order to obtain operations starting and finishing in the twistor spaces themselves,
we just need to recall \cite[Section 7]{LB91} that the twistor space can be obtained from the projective model by (a) blowing-down the two distinguished sections $E$ and $\ol E$, and (b) taking small resolutions for all double points.
That is to say, it is enough to attach these operations beforehand and afterward.
However, once we can identify the output space $X'$ (as we did in Theorem \ref{thm-3}), 
it is possible to restate our construction without passing to projective models as follows:

\begin{enumerate}
\item[1.] Let $Z$ be any LeBrun twistor space, and $\Phi:Z\to \mathbb{CP}^3$  the rational map associated to the fundamental system $|K^{-1/2}|$.
The image $\Phi(Z)$ is a non-singular quadric $Q$, and 
for any discriminant curve $\mathscr C_i=\{P_i=0\}\subset Q$, the inverse image $\Phi^{-1}(\mathscr C_i)$ consists of two irreducible components $D_i$ and $\ol D_i$.

\item[2.] As before we move any one of the discriminant curves, say $\mathscr C_1$, to a reducible curve $\mathscr C_1^{\infty}=\{uv=0\}$.
(In order to justify this process rigorously, we need the projective models as we discussed.)
Let $Z'$ be the limit variety.
Then in $Z'$ each of $D_i$ and $\ol D_i$ breaks into two irreducible components.
Among these two components, exactly one is biholomorphic to $\mathbb{CP}^2$.
So over the reducible curve $\mathscr C_1^{\infty}$ we have two $\mathbb{CP}^2$-s.
The intersection of these $\mathbb{CP}^2$-s is a (real) ordinary double point of $Z'$.

\item[3.] We take the small resolution of this ordinary double point which does not change the structure of the two $\mathbb{CP}^2$-s. (This condition uniquely specifies the small resolution.)
Then the normal bundle of these becomes isomorphic to $\mathscr O(-1)$.

\item[4.] So we blow-down the two $\mathbb{CP}^2$-s.
Then the resulting space is nothing but the LeBrun twistor space on $(n-1)\mathbb{CP}^2$, whose projective model is exactly the variety $X'$ in Theorem \ref{thm-3}.
\end{enumerate} 

A proof of these can be readily obtained if we notice that 
in each steps in the former construction that uses the projective models,
the two divisors $E$ and $\ol E$ are always contained in such a way that 
they can be blown-down to $\mathbb{CP}^1$,
as one can read off from Figures \ref{fig-n_1}--\ref{fig-n_3}.
So we omit the proof.

\begin{remark}
{\em
The result in this section gives a way to obtain a LeBrun twistor on $(n-1) \mathbb{CP}^2$ from that on $n \mathbb{CP}^2$.
This construction can be readily generalized to give  an explicit way for obtaining a LeBrun twistor on $(n-k)\mathbb{CP}^2$ from that on $n \mathbb{CP}^2$, for any $1\le k\le n$.
Namely, among $n$ discriminant curves $\mathscr C_1,\cdots,\mathscr C_n$, we choose any $k$ curves, say $\mathscr C_1,\cdots,\mathscr C_k$.
Then we  move each of these curves $\mathscr C_i$ to reducible curves $\mathscr C_i^{\infty}$ as in the case $k=1$.
We suppose that all intersection points of the irreducible components of these reducible curves are transversal.
The inverse images of these $(1,0)$- and $(0,1)$-curves consists of two irreducible components,
which are isomorphic to $\Sigma_1$ and $\Sigma_{n-1}$.
So in all we have $2k$ $\Sigma_1$-s and $2k$ $\Sigma_{n-1}$-s.
They intersect in such a way that any two different $\Sigma_1$-s  shares a unique point, and the same for any two different $\Sigma_{n-1}$-s.
(These are ordinary double points of the limit 3-fold.)
Hence we can choose a small resolution of these ordinary double points which does not change $\Sigma_1$.
Consequently we are again in the situation that these $2k$ $\Sigma_1$-s can be blown-down to $\mathbb{CP}^1$,
which gives a projective model of LeBrun twistor space on $(n-k)\mathbb{CP}^2$. 

}
\end{remark}
\section{A degeneration to LeBrun metrics on $\mathscr O(-n)$ }
In \cite[Section 5]{LB91} LeBrun showed by direct calculations for metrics that if the monopole points $p_1,\cdots,p_n$ on $\mathscr H^3$ converge to a single point, then
after a suitable conformal change, the metrics  converge to the scalar-flat K\"ahler metric on the total space of the line bundle
$\mathscr O(-n)\to\mathbb{CP}^1$ constructed also by LeBrun  \cite{LB88}.
In this section we prove that if all $n$ polynomials $P_1,\cdots,P_n$ (corresponding to the monopole points) are equal, 
then the variety obtained from the space $X$  defined in \eqref{LB-n} by blowing-down  the divisors $E$ and $\ol{E}$
 is (a compactification of) the  twistor space of the LeBrun  metric on $\mathscr O(-n)$.
  We show this by first finding all twistor lines explicitly, and then use  a result by Pontecorvo \cite{Pont} and a characterization of the LeBrun metric showed in \cite{LB88}.

In order to reduce complexity in computation, we use slightly different
coordinates on  $\mathbb{CP}^1\times\mathbb{CP}^1$ from  Section 2.
Namely we adapt non-homogenous coordinates  $(u,v)$ on  $\mathbb{CP}^1\times\mathbb{CP}^1$ with respect to which the real structure is given by $(u,v)\mapsto (-1/\ol v,-1/\ol u)$,
and the unique polynomial becomes $P_1=u-v$.
So the variety $X$ in \eqref{LB-n} becomes
\begin{equation}\label{LB-n''}
X_0:=\{(x,y,z)\in \mathbb P(\mathscr E_n)
\set
xy=(u-v)^nz^2
\}.
\end{equation}
This is invariant under the real structure on $\mathscr E_n$ defined by
\begin{equation}\label{rsq''}
(x,y,z)\longmapsto \
\left(
\frac{\ol{y}}{\ol{u}^{n-1}\ol{v}},\,
(-1)^n\frac{\ol{x}}{\ol{u}\,\ol{v}^{n-1}},\, \ol{z}\right),
\end{equation}
which is a lift of the real structure on $\mathbb{CP}^1\times\mathbb{CP}^1$.
It is elementary to see that under the real structure \eqref{rsq''} the variety $X_0$ has no real point.

\begin{remark}{\em
When $n=2$, the variety $X_0$ is isomorphic to the variety $X$ in \eqref{LB-2} with
$\lambda=1$.
Recalling that $\lambda\in(0,1)$ represents a self-dual conformal structure  of positive scalar curvature on $2 \mathbb{CP}^2$,
$X_0$ corresponds to another limit in the moduli space. 
}
\end{remark}

Of course $X_0$ still has a conic bundle structure $f:X_0\to Q$ and its discriminant locus is the curve $\Delta:=\{u=v\}$.
The inverse image $f^{-1}(\Delta)$ consists of two components $\{x=0\}$ and $\{y=0\}$.
We further  define
\begin{equation}
L_{\infty}:=\{x=y=u-v=0\}
\end{equation}
which is precisely the intersection of the two components.
$L_{\infty}$ is a real smooth rational curve.
$X_0$ has $A_{n-1}$-singularities along $L_{\infty}$ and is non-singular outside $L_{\infty}$.
$X_0$ still has two section $E=\{x=z=0\}$ and $\ol{E}=\{y=z=0\}$ and they can be blow-down
to $\mathbb{CP}^1$ by the same reasoning as before.
Let $\mu:X_0\to Z_0$ be the blowing-down.
Since $L_{\infty}$ is disjoint from $E$ and $\ol{E}$, $L_{\infty}$ is naturally contained in $Z_0$.
We denote this by the same notation $L_{\infty}$.
Then we have the following

\begin{theorem}\label{thm-4}
The complement $Z_0\backslash L_{\infty}$ has a structure of the twistor space of a scalar flat K\"ahler metric on the total space of $\,\mathscr O(-n)$.
\end{theorem}
\begin{proof}
We prove the claim by finding all twistor lines in explicit form.
It is elementary to deduce that any smooth $(1,1)$-curve  contained in $Q\backslash Q^{\sigma}$ which is invariant under the real structure on $Q$ must be of the form, in the above coordinate $(u,v)$, 
\begin{equation}\label{mtl}
u(t)=\frac{d-rt}{1+r\ol{d}t},\hspace{2mm}
v(t)=\frac{rd-t}{r+\ol{d}t},
\end{equation}
where $t\in\mathbb C\cup\{\infty\}$ is a coordinate on the curve, $d$ and $r$ are constants with $d\in\mathbb C\cup\{\infty\}$ and $0<r<1$
When $d=\infty$, \ \eqref{mtl} means $u=1/(rt),\,v=r/t$.
We denote the curve \eqref{mtl} by $\mathscr C(d,r)$.
When $d\neq \infty$, from \eqref{mtl} we compute
\begin{equation}
u(t)-v(t)=\frac{(1+|d|^2)(1-r^2)t}{(1+r\ol{d}t)(r+\ol{d}t)}.
\end{equation}
Therefore the restriction of the conic bundle $X_0\to Q$ onto $\mathscr C(d,r)$ is explicitly given by the equation
\begin{equation}\label{rest1}
\xi\eta=\frac{(1+|d|^2)^n(1-r^2)^nt^n}{(1+r\ol{d}t)^n(r+\ol{d}t)^n},
\end{equation}
where, we are working on the affine set $z\neq 0$ and writing $\xi=x/z\in\mathscr O(n)$ and $\eta=y/z\in\mathscr O(n)$.
Viewing as an equation for $(\xi,\eta)$ with $\xi$ and $\eta$ being rational functions of $t$, \eqref{rest1} has  solutions 
\begin{equation}\label{sol1}
\xi(t)=c\,\frac{(1+|d|^2)^{\frac n2}(1-r^2)^{\frac n2}t^n}{(1+r\ol{d}t)(r+\ol{d}t)^{n-1}},
\,\,\eta(t)=c^{-1}\frac{(1+|d|^2)^{\frac n2}(1-r^2)^{\frac n2}}{(1+r\ol{d}t)^{n-1}(r+\ol{d}t)},
\end{equation}
where $c\in {\rm{U(1)}}$ is arbitrary constant.
Note that since  $0<r<1$, the square root of $1-r^2$ has a unique meaning;
namely we choose a positive one.
Let $L(d,r,c)$ be the real curve  \eqref{sol1}.
Thanks to  $c\in {\rm{U(1)}}$,  this curve is invariant under the real structure \eqref{rsq''}.
Similarly, when $d=\infty$, if we define the curve $L(\infty, r,c)$ over $\mathscr C(\infty,r)$ by
\begin{equation}
\xi(t)=\frac{(1-r^2)^{\frac n2}}{r},
\quad
\eta(t)=\frac{(1-r^2)^{\frac n2}}{r^{n-1}t^n},
\end{equation}
then $L(\infty,r,c)$  is real and
contained in $X_0$.
From the explicit forms, all these curves $L(d,r,c)$, $d\in\mathbb C\cup\{\infty\}, \,0<r<1,\,c\in{\rm{U(1)}}$, are disjoint from the curve $L_{\infty}$.
Thus we  obtain a 4-dimensional family 
\begin{equation}
\mathscr L:=\left\{L(d,r,c)\set d\in\mathbb C\cup\{\infty\},\, 0<r<1,\, c\in {\rm U(1)}\right\}
\end{equation}
of real curves in $X_0$.

We now prove the claim that {\em the complement
\begin{equation}\label{compl}
X_0\backslash(E\cup\ol{E}\cup L_{\infty}\cup f^{-1}(Q^{\sigma}))
\end{equation}
is foliated by members of $\mathscr L$}.
Namely we show that for any point of this set there exists a unique member of  $\mathscr L$ going through the point.
We prove this by taking arbitrary point of $Q\backslash Q^{\sigma}$  and verify that the claim holds for any  point of the fiber over there, as long as the point is  not  on $E\cup\ol{E}\cup L_{\infty}$.

We begin with the easiest case. 
Take $(0,0)\in\Delta$.
Its fiber is $\{xy=0\}$, a union of two lines.
It can be obtained by direct computations that the curve $\mathscr C(d,r)$ goes through $(0,0)$ iff $d=0$ or $d=\infty$.
Similarly using \eqref{sol1}  we deduce that the intersection point of $L(0,r,c)$ with the fiber is given by
\begin{equation}
\xi=0,\,\,\eta=\frac{c^{-1}(1-r^2)^{\frac n2}}{r}.
\end{equation}
Similarly the intersection point of $L(\infty,r,c)$ with the fiber is given by
\begin{equation}
\xi=\frac{c(1-r^2)^{\frac n2}}{r},\,\,\eta=0.
\end{equation}
Further it can be readily shown that the function $r\mapsto (1-r^2)^{n/2}/r$ is 
strictly decreasing on the interval $(0,1)$, tends to $+\infty$ when $r\to 0$, and tends to $0$ when $r\to 1$.
This implies that any point of the fiber, except the 3 points $(x:y:z)=(0:1:0)\in E,\, (1:0:0)\in \ol{E}$ and $(0:0:1)\in L_{\infty}$, is passed by a unique member of $\mathscr L$ (of the form $L(0,r,c)$ if $\xi=0$ and $L(\infty,r,c)$ if $\eta=0$).
Thus we obtain the claim, for the fiber over $(0,0)\in Q$.

Next we choose a point $(0,v)\in Q\backslash(Q^{\sigma}\cup\Delta)$.
This assumption implies $v\in\mathbb C^*$. 
The fiber over $(0,v)$ is a smooth conic $xy=(-v)^nz^2$.
We show that for any point of the conic satisfying $z\neq 0$, there is a unique member of $\mathscr L$ going through the point.
Suppose $(0,v)\in\mathscr C(d,r)$ for some $(d,r)\in(\mathbb C\cup\{\infty\})\times (0,1)$.
Then we readily obtain $d\not\infty$ and, from \eqref{mtl} we get
\begin{equation}
v(d/r)=\frac{d(r^2-1)}{r^2+|d|^2}=v.
\end{equation}
Solving this with respect to $r^2$ we obtain
$
r^2=d(1+\ol{d}v)/(d-v).
$
Since this is a real number, we obtain
\begin{equation}
\frac{d(1+\ol{d}v)}{d-v}=\frac{\ol{d}(1+d\ol{v})}{\ol{d}-\ol{v}}.
\end{equation}
This is equivalent to $(1+|d|^2)(\ol{d}v-d\ol{v})=0$ and hence we obtain $d\ol{v}\in\mathbb R$.
So we put $d\ol{v}=a\in\mathbb R$. 
Then we have
\begin{equation}
r^2=\frac{a(1+a)}{a-|v|^2}.
\end{equation}
From this we obtain that being 
 $r\in (0,1)$ is equivalent to $a\in (-1,0)$.
Hence we have obtained that  the curve $\mathscr C(d,r)$ goes through the point $(0,v)$ 
(with $v\in\mathbb C^*$) iff there is a constant $-1<a<0$ such that 
\begin{equation}\label{dr}
d=\frac{a}{\ol{v}},\,\,r=\sqrt{\frac{a(1+a)}{a-|v|^2}}
\end{equation}
hold.
Further, the constant $-1<a<0$ is uniquely determined by $(d,r)$.
Substituting \eqref{dr} into \eqref{sol1} and considering the absolute value, we obtain
\begin{equation}
|\xi|^2=\frac{-a(|v|^2-a)^{n-1}}{1+a}.
\end{equation}
Viewing this as a function of $-1<a<0$, it is readily seen that it tends to $+\infty$ when $a\to  -1$ and goes to zero if $a\to 0$, and its differential is always negative.
This implies that for arbitrary $\xi\in\mathbb C^*$ there exists a unique element $(d,r,c)\in \mathbb C\times (0,1)\times U(1)$ such that the $\xi$-coordinate of  $L(d,r,c)$ coincides with $\xi$.
Since a point of the fiber is uniquely determined by the value of the $\xi$-coordinate, 
it follows  that any point of the fiber over $(0,v)$ is passed by a unique member of $\mathscr L$, except (possibly) on the 2 points $(1,0,0)\in E$ and $(0,1,0)\in \ol E$.

As the final step of the verification of the claim, 
we make use of automorphisms of $X_0$.
First, we note that holomorphic automorphisms of $Q$ commuting with the real structure   and preserving the discriminant locus $\Delta$ must be of the form 
\begin{equation}\label{su2}
(u,v)\longmapsto
\left(
\frac{\alpha u+\beta}{-\ol{\beta}u+\ol{\alpha}},
\frac{\alpha v+\beta}{-\ol{\beta}v+\ol{\alpha}}
\right),
\end{equation}
where $\alpha,\beta\in\mathbb C$ with $|\alpha|^2+|\beta|^2=1$,
which constitute an SU(2)-action.
It can be shown by direct computations that this action  lifts on the line bundles $\mathscr O(n-1,1)$, $\mathscr O(1,n-1)$ and $\mathscr O$ in such a way that
it preserves the variety $X_0$, and it commutes with the real structure \eqref{rsq''}.
(We do not write the explicit form.)
The SU(2)-action  on $X_0$ thus obtained maps members of $\mathscr L$ to those of $\mathscr L$.

Take arbitrary point $(u,v)$ of $Q\backslash Q^{\sigma}$.
Then by \eqref{su2} there exists a real automorphism $\phi$ keeping $\Delta$ such that  $\phi(u,v)=(0,v')$ for some $v'\in\mathbb C$.
These implies that
if some two different members of $\mathscr L$ intersect on the fiber over  $(u,v)$, then there must be a pair of members of $\mathscr L$ intersecting on the fiber over $(0,v')$.
This contradicts what we have verified for the fibers over the point $(0,v)$, $v\in\mathbb C$.
By similar reasoning (using the automorphism \eqref{su2}), we obtain that 
for any point on the fiber over arbitrary point $(u,v)\in Q\,\backslash Q^{\sigma}$, 
there exists a member of $\mathscr L$ going through the point.
Thus we have obtained that the family $\mathscr L$ actually foliates the set \eqref{compl}.
In the following $\mathscr L$ means the set of image curves  under the blowing-down  $\mu:X_0\to Z_0$.

On the other hand, for any point of $Q^{\sigma}$ its fiber is obviously real.
(Recall that $\sigma$ has no real point.)
Let $\mathscr F$ be the set of these real fibers of $f$.
Taking the images under $\mu$,
we obtain a set of real curves in $Z_0$ parameterized by $Q^{\sigma}$.
We still write $\mathscr F$ for this family.

Thus we have obtained two families $\mathscr L$ and $\mathscr F$ of  real smooth rational curves in $Z_0$.
From the construction the union $\mathscr L\cup\mathscr F$ foliates $Z_0\backslash L_{\infty}$, where
 any point on the curves $\mu(E)$ and $\mu(\ol{E})$ is passed by a unique member of $\mathscr F$.
In order to show that $Z_0\backslash L_{\infty}$ actually has a structure of a twistor space it remains to see that the normal bundles of these curves are all isomorphic to $\mathscr O(1)^{\oplus 2}$. 
This is immediate to see for members of $\mathscr F$.
On the other hand any member of $\mathscr L$ is an image of some $L=L(d,r,c)$ in $X_0$.
By construction $L\subset f^{-1}(\mathscr C(d,r))$ holds.
Further it is easily seen that the surface  $ f^{-1}(\mathscr C(d,r))$ has  $A_{n-1}$-singularities over the intersection points $\Delta\cap\mathscr C(d,r)$ and is smooth outside these 2 points.
In particular, as $L\cap L_{\infty}=\emptyset$, the surface is smooth in a neighborhood of $L$.
Then since $L$ can actually be moved (by changing $c\in U(1)$) in such a way that it is disjoint from $L$,  the normal bundle of $L$ in $f^{-1}(\mathscr C(d,r))$ is trivial.
It follows from the two inclusions $L\subset f^{-1}(\mathscr C(d,r))\subset X_0$ that the normal bundle of $L$ in $X_0$ is either $\mathscr O(1)^{\oplus 2}$ or $\mathscr O(2)\oplus \mathscr O$.

We show that  the latter cannot happen.
As the normal bundle is $\mathscr O(1)^{\oplus 2}$ or $\mathscr O(2)\oplus\mathscr O$,
whose $H^1$-s are zero,
the moduli space of all small deformations of $L$ in $X_0$ is non-singular complex 4-manifold and its tangent space at the point (corresponding to $L$) is naturally isomorphic to $H^0(N)$, where $N$ is the normal bundle.
Because  $\mathscr L$ is a real 4-dimensional family, 
it must be a real slice of the moduli space.
If $N$ is isomorphic to $\mathscr O(2)\oplus\mathscr O$, then 
there is a real section $s\in H^0(N)$ such that its zeros consist of distinct 2 points.
Further,  since the surface $f^{-1}(\mathscr C(d,r))$ is real in $X_0$, 
the section $s\in H^0(N)$ can be supposed to be real.
Let $p$ and $\ol p\in L$ be the zeros of such an $s$.
Then $s$ generates a complex 1-dimensional deformation of $L$ in $X_0$ such that 
all members go through $p$ and $\ol p$.
Further, by reality of $s$, this  deformation contains a real 1-dimensional subfamily whose members are real.
This contradicts the fact that any two different members of the family $\mathscr L$ are disjoint.
Hence we obtain $N\simeq\mathscr O(1)^{\oplus 2}$.
Then since $\mu$ is isomorphic outside $E\cup \ol{E}$, the normal bundle of $\mu(L)$ in $Z_0$ remains to be $\mathscr O(1)^{\oplus 2}$.
Thus we obtain that $Z_0\backslash L_{\infty}$ is foliated by real smooth rational curves without real point whose normal bundles are isomorphic to $\mathscr O(1)^{\oplus 2}$.

Next we see that the parameter space of these twistor lines is diffeomorphic to $\mathscr O(-n)$. 
As is already mentioned the inverse image $f^{-1}(\Delta)$ consists of two divisors $\{x=0\}$ and $\{y=0\}$. 
These are clearly biholomorphic to a ruled surface $\Sigma_n$, and  their intersection $L_{\infty}$ is a section of the ruling whose self-intersection  is $(+n)$.
Thus both of the connected component of $f^{-1}(\Delta)\backslash L_{\infty}$ is biholomorphic to the total space of $\mathscr O(-n)$.
Further it follows from our concrete description that the twistor lines we have detected intersect   each of the connected components of $f^{-1}(\Delta)\backslash L_{\infty}$ at a unique point.
Hence the parameter space of the twistor lines, namely the base 4-manifold, has to be diffeomorphic to $\mathscr O(-n)$.

Finally we have to show that the self-dual conformal structure on $\mathscr O(-n)$ is represented by a scalar-flat K\"ahler metric (after reversing the orientation).
By adjunction formula and the inclusion $X_0\subset \mathbb P(\mathscr E_n)$ we readily obtain, outside the singular locus $L_{\infty}$, that 
$K_{X_0}\simeq
\mathscr O_{X_0}(-E-\ol E)\otimes
f^* \mathscr O(-2,-2)$.
From this we obtain 
$\mu^*K_{Z_0}^{-1/2}
\simeq 
\mathscr O_{X_0}(E+\ol E)
\otimes 
f^*\mathscr O(1,1)$
on $X_0\backslash L_{\infty}$.
As $\Delta$ is a $(1,1)$-curve, this means
$\mu(f^{-1}(\Delta)\backslash L_{\infty})\in |K_{Z_0\backslash L_{\infty}}^{-1/2}|$.
Therefore by a theorem of Pontecorvo \cite{Pont}, 
the conformal class is represented by a scalar-flat K\"ahler metric (on $\mathscr O(-n)$).
Thus we have completed a proof of Theorem \ref{thm-4}.
\end{proof}

With the aid of the characterization result by LeBrun \cite{LB88} we obtain the following

\begin{theorem}\label{thm:O(-n)}
$Z_0\backslash L_{\infty}$ is the twistor space of LeBrun's asymptotically flat, scalar flat K\"ahler  metric on\, $\mathscr O(-n)$ constructed in \cite{LB88}.
\end{theorem}	

\begin{proof}
By pulling back the metric obtained in Theorem \ref{thm-4}
under the usual $n$-fold covering $\mathscr O(-1)\to\mathscr O(-n)$, we get a scalar-flat K\"ahler metric on $\mathscr O(-1)\backslash\{0\}\,(\simeq\mathbb C^2\backslash\{0\})$.
Let $W$ be the twistor space of this metric,
and $\pi:W\to Z_{0}\backslash L_{\infty}$  the natural $n$-fold covering.
Since $Z_0$ is locally a product of a surface $A_{n-1}$-singularity and a disk ($\subset \mathbb C)$ in a neighborhood of each point of $L_{\infty}$,
and the surface $A_{n-1}$-singularity is uniformized by an $n$-fold covering,
we can add a curve $\tilde L_{\infty}\simeq L_{\infty}$ to $W$ to make
the union $\ol W:=W\cup \tilde L_{\infty}$ non-singular 3-fold
and the projection $\ol{\pi}:\ol W\to Z_0$ give a uniformization
for a neighborhood of $L_{\infty}$ in $Z_0$.

Let $D$ and $\ol D$ be the two irreducible components of $\mu(f^{-1}(\Delta))$.
Then as is evident from the defining equation of $D$ and $\ol D$, 
the inverse images $\ol{\pi}^{-1}(D)$ and $\ol{\pi}^{-1}(\ol D)$ are also non-singular surfaces
 containing $\tilde L_{\infty}$ and intersecting transversally along $\tilde L_{\infty}$.
 Further, since
 $N_{L_{\infty}/D}\simeq N_{L_{\infty}/\ol D}\simeq \mathscr O(n)$ and 
  $D$ and $\ol D$ are precisely the image of $\ol{\pi}^{-1}(D)$ and $\ol{\pi}^{-1}(\ol D)$ 
 under the quotient map by $\mathbb Z_n$-action which fixes $\tilde L_{\infty}$, the normal bundles of $\tilde L_{\infty}$ in  $\ol{\pi}^{-1}(D)$ and $\ol{\pi}^{-1}(\ol D)$ have to be $(+1)$.
Therefore, the normal bundles of $\tilde L_{\infty}$ in $\ol W$ are isomorphic to $\mathscr O(1)^{\oplus 2}$.
Therefore the 3-fold $\ol W$ is a twistor space of a 4-manifold $(\mathscr O(-1)\backslash\{0\})\cup\{\tilde{\infty}\}$,
where $\tilde{\infty}$ is the point corresponding to $\tilde L_{\infty}$.

Since the divisor $D+\ol D\in |K_{Z_0}^{-1/2}|$ outside $L_{\infty}$,
we have $\ol{\pi}^{-1}(D)+\ol{\pi}^{-1}(\ol D)\in |K_{W}^{-1/2}|$.
Therefore by the proof of Proposition 6 in \cite{LB92}, 
the conformal class on  $\mathscr O(-1)\backslash \{0\}$ determined by the twistor space $W$ is represented by an {\em asymptotically flat}, scalar flat K\"ahler metric.
Hence we have shown that the scalar-flat K\"ahler metric in Theorem \ref{thm-4} is 
locally asymptotically flat.
 On the other hand  as seen in the proof of Theorem \ref{thm-4}, 
 the space $W$ admits an effective SU(2)-action commuting with the real structure.
 Hence the metric admits an effective SU(2)-action.
 Therefore by characterization result of LeBrun
\cite[p.595, Proposition]{LB88} the metric must be LeBrun's metric, as desired.
\end{proof}

\section{A degeneration to the Gibbons-Hawking metrics}

In  \cite{GH78}, Gibbons-Hawking   explicitly constructed a family of hyperK\"ahler metrics on  complex surfaces which are diffeomorphic to the minimal resolution of $\mathbb C^2/\Gamma$, where $\Gamma$ is a cyclic subgroup of SU(2) of any finite order.
In \cite{Hi79}, Hitchin constructed the twistor spaces of these metrics.
In this section we show that these twistor spaces can be obtained as a limit of a Zariski open subset of LeBrun's twistor spaces on $n \mathbb{CP}^2$, where $n$ is exactly the order of the cyclic subgroup $\Gamma$.

We begin by quickly recalling the Hitchin's construction.
%
Let $p_i(u)=a_iu^2+2b_iu+c_i$, $1\le i\le n$, be  distinct $n$ quadratic polynomials of $u$, with the coefficients satisfying  reality conditions 
$a_i=-\ol c_i,{\text{ and }} b_i\in \mathbb R.
$
Denote by $\mathscr O(d)$  the line bundle of degree $d$ over $\mathbb{CP}^1$.
Let $\tilde Z$ be  a 3-dimensional algebraic variety defined by
\begin{align}\label{Htw}
\tilde Z:=\left\{(x,y,z)\in \mathscr O(n)^{\oplus 2}\oplus \mathscr O(2)\set xy=(z-p_1(u))\cdots(z-p_n(u))\right\},
\end{align}
where $(x,y)\in\mathscr O(n)^{\oplus 2}$, $z\in \mathscr O(2)$. 
By using the projection to $\mathscr O(2)$, we regard the variety $\tilde Z$ as    a  conic bundle over the total space of $\mathscr O(2)$, 
where the conics are in $\mathbb C^2$.
Then the discriminant curves of the conic bundle are exactly the sections $\{z=p_i(u)\set 1\le i\le n\}$. 
Therefore, the variety $\tilde Z$ has singularity precisely over the intersection points of these sections.
Moreover, by the reality condition for the coefficients of $p_i(u)$,  
any two different sections intersect transversally at (mutually conjugate) two points.
This means that all the singularities of $\tilde Z$ admit a small resolution.
Furthermore, $\tilde Z$ has a real structure naturally induced from the anti-podal map on $\mathbb{CP}^1$. 
Then the desired twistor space of the Gibbons-Hawking metric is obtained as an appropriate small resolution of $\tilde Z$, equipped with a natural lift of the real structure.

For giving a precise relationship between these twistor spaces with LeBrun's ones,
it is convenient to think the total space of the above line bundle $\mathscr O(2)$ as a degeneration of
$Q=\mathbb{CP}^1\times\mathbb{CP}^1$, which is the minitwistor space  of the hyperbolic space $\mathscr H^3$.
More explicitly, we define a subvariety $\mathscr Q\subset\mathbb{CP}^3\times \mathbb C$ by
\begin{align}
\mathscr Q:=\left\{(x_0,x_1,x_2,x_3)\times s\in\mathbb{CP}^3\times \mathbb C
\set
s^2x_0^2-x_1^2-x_2x_3=0\right\}.
\end{align}
We also define a real structure on $\mathscr Q$ by
\begin{align}\label{rsQ}
(x_0,x_1,x_2,x_3)\times s
\longmapsto
(\ol x_0,\ol x_1, \ol x_3,\ol x_2)\times \ol s.
\end{align}
Let $\pi:\mathscr Q\to\mathbb C$ be the projection to the last factor.
If we view this as a family of quadrics in $\mathbb{CP}^3$ parameterized by $s\in\mathbb C$, 
then a fiber $Q_s:=\pi^{-1}(s)$ is non-singular when $s\neq 0$, and a cone over a conic $\{x_1^2+x_2x_3=0\}\subset\mathbb{CP}^2$ 
with vertex ${\bm o}:=(1,0,0,0)$ when $s=0$.
Note that the variety $\mathscr Q$ has a unique singularity at ${\bm o}$, and it is an ordinary double point.
Also, a fiber $Q_s$ is invariant under the real structure \eqref{rsQ} if and only if $s\in\mathbb R$.
In particular, the cone $Q_0$ itself is real, and we remark that on the conic $\{x_1^2+x_2x_3=0\}$, the real structure is naturally acting as an anti-podal map since there is no real point on the conic.
Therefore, the real locus on $Q_0$ is the vertex only.
On the other hand, if $s\in\mathbb R^{\times}$, then the real locus on $Q_s$ is a (non-holomorphic) 2-sphere.
If $s$ tends to zero, this sphere shrinks to the vertex.

As is well-known, the family  $\pi:\mathscr Q\to\mathbb C$ can be trivialized outside the origin.
Actually, by rewriting the defining equation of $\mathscr Q$  as
$$
(sx_0-x_1)(sx_0+x_1)=x_2x_3,
$$
we can take the trivializing map $\phi:\mathscr Q|_{\pi^{-1}(\mathbb C^*)}\simeq (\mathbb{CP}^1\times\mathbb{CP}^1)\times\mathbb C^*$ by
%
\begin{align}\label{isom15}
(x_0,x_1,x_2,x_3)\times s&\stackrel{\phi}{\longmapsto}&(sx_0-x_1,x_2)\times(sx_0-x_1,x_3)\times s\\
&=&(x_3,sx_0+x_1)\times
(x_2,sx_0+x_1)\times s.
\end{align}

Next we give two (non-Cartier) divisors on the variety $\mathscr Q$ by
\begin{align}\label{scrD1}
\mathscr D_1:=\left\{(x_0,x_1,x_2,x_3)\times s\in\mathbb{CP}^3\times \mathbb C
\set
x_3=sx_0-x_1=0
\right\},
\end{align}
\begin{align}\label{scrD2}
\mathscr D_2:=\left\{(x_0,x_1,x_2,x_3)\times s\in\mathbb{CP}^3\times \mathbb C
\set
x_2=sx_0-x_1=0
\right\}.
\end{align}
Then clearly these are contained in $\mathscr Q$, and
for any $s\in\mathbb C$, the intersections $Q_s\cap \mathscr D_1$ and $Q_s\cap \mathscr D_2$
are lines in $\mathbb{CP}^3$.
Further, if $s\neq 0$, the restrictions $\mathscr D_1|_{Q_s}$ and $ \mathscr D_2|_{Q_s}$ have mutually different bidegree  $(1,0)$ and $(0,1)$ for the identification $Q_s\simeq\mathbb{CP}^1\times\mathbb{CP}^1$ given by  \eqref{isom15}, and  if $s=0$ both  become an identical line (an edge of the cone).

While $\mathscr D_1$ and $\mathscr D_2$ are non-Cartier divisors on $\mathscr Q$, 
we can consider the associated line bundles $[\mathscr D_1]$ and $[\mathscr D_2]$, at least over the non-singular locus $\mathscr Q\backslash\{{\bm o}\}$.
Then the following property of these line bundles is obvious from the above consideration:

\begin{lemma}\label{lemma:rest4}
Under the trivializing map $\phi$ over $\mathbb C^*$ for the family $\mathscr Q$
 (see \eqref{isom15}), for any $s\neq 0$, we have the following isomorphisms of line bundles:
$$[\mathscr D_1]|_{Q_s}\simeq \mathscr O(1,0),
\quad
[\mathscr D_2]|_{Q_s}\simeq \mathscr O(0,1).
$$
Further, for $s=0$, letting $\bm e$ be any edge of the cone $Q_0$, we have
$$
[\mathscr D_1]|_{Q_0\backslash\{{{\bm o}}\}}\simeq
[\mathscr D_2]|_{Q_0\backslash\{{{\bm o}}\}}\simeq \mathscr O_{Q_0\backslash\{{{\bm o}}\}}([\bm e]).
$$
\end{lemma}

With these preliminaries, we take and {\em fix} any LeBrun twistor space on $n \mathbb{CP}^2$.
As explained in Section \ref{s:2}, this is equivalent to taking $n$ discriminant curves $\mathscr C_1,\cdots,\mathscr C_n$ on $Q=\mathbb{CP}^1\times\mathbb{CP}^1$, which are of bidegree $(1,1)$ satisfying the relevant properties.
{\em If we fix an embedding of $Q$ into $\mathbb{CP}^3$} as a quadric,  any such a curve is written as a unique plane section $Q\cap H$.


\begin{lemma}\label{lemma:hi}
If we realize $Q$ as $Q_s\subset\mathbb{CP}^3$ $(s\neq0)$,
the plane $H$ corresponding to a point on $\mathscr H^3$ is defined 
by the equation of the form
\begin{align}\label{H_i}
x_0-bx_1-cx_2-\ol c x_3=0,
\end{align}
where $b\in \mathbb R$ and $c\in\mathbb C$ satisfying the inequality
\begin{align}\label{ellipsoid}
b^2+4|c|^2<1/s^2.
\end{align}
In particular, the plane cannot not hit the vertex $\bm o=(1,0,0,0)$.
\end{lemma}

Note that the inequality \eqref{ellipsoid} involves the parameter $s$.
As this is deduced by direct computation in a similar way to \cite[Lemma 6.11]{HV09}, we omit a proof and just remark that the constraint \eqref{ellipsoid} is a consequence of the absence of real point on $H\cap Q_s$.

For the chosen LeBrun metric, we take an embedding $Q=Q_1$ (i.e.\,the quadric $Q_s$ in the case $s=1$),
and  let $h_1, \cdots,h_n$ be the defining linear polynomials of $H_1,\cdots,H_n$ respectively.
Each $h_i$ can be regarded as a section of $\mathscr O_{\mathbb{CP}^3}(1)$,
so that by pulling back and restricting to the subvariety $\mathscr Q$, we obtain 
$n$ sections of 
$\mathscr O_{\mathscr Q}(1)
:=
(p^*\mathscr O_{\mathbb{CP}^3}(1))|_{\mathscr Q}$,
where $p:\mathbb{CP}^3\times\mathbb C\to \mathbb{CP}^3$
is the projection to the first factor.
Furthermore, for this line bundle, since in $\mathscr Q$ the sum $\mathscr D_1+\mathscr D_2$ is defined by a single linear equation $sx_0-x_1=0$,  outside the vertex $\bm o$, we have  an isomorphism
\begin{align}\label{isom111}
\mathscr O_{\mathscr Q}(1)
\simeq 
[\mathscr D_1]\cdot[\mathscr D_2]
\quad
\text{(tensor product)}.
\end{align}
Therefore we have
\begin{align}\label{rest111}h_i|_{\mathscr Q\backslash\{\bm o\}}\in H^0([\mathscr D_1]\cdot[\mathscr D_2]).
\end{align}
We also remark that the restriction $h_i|_{Q_s}$ $(s\neq 0$) makes sense
and via the trivializing map $\phi$ it can be regarded as a section of $\mathscr O(1,1)$.

Next we define two line bundles over $\mathscr Q\backslash\{\bm o\}$ by
$\mathscr L_1:=[\mathscr D_1]^{n-1}\cdot[\mathscr D_2]$ and 
 $\mathscr L_2:=[\mathscr D_1]\cdot[\mathscr D_2]^{n-1}$.
Then we shall define 
\begin{align}\label{eqn:scrZ}
\mathscr Z:=\left\{(x,y)\in 
\mathscr L_1
\oplus
\mathscr L_2
\set
 xy=h_1h_2\cdots h_n\right\},
\end{align}
and let $q:\mathscr Z\to\mathscr Q\backslash\{\bm o\}$ be the projection.
Then by $q$, the space $\mathscr Z$ has a structure of a conic bundle over $\mathscr Q\backslash\{\bm o\}$, where conic again means a quadratic curve in  $\mathbb C^2$.
Here we note that the (tensor) product $xy$ belongs to 
the line bundle $[\mathscr D_1]^n\cdot[\mathscr D_2]^n$, which is isomorphic to 
$\mathscr O_{\mathscr Q}(n)$ by \eqref{isom111}.
On the other hand, by \eqref{rest111}
the product $h_1h_2\cdots h_n$ belongs to $H^0(\mathscr O_{\mathscr Q}(n))
\simeq 
H^0([\mathscr D_1]^n\cdot[\mathscr D_2]^n)$.
These mean that the locus $\{h_i=0\set 1\le i\le n\}$ are all the discriminant locus of the conic bundle $q$.

For the real structure, by noting that the real structure \eqref{rsQ} interchanges the divisors $\mathscr D_1$ and $\mathscr D_2$, we have a natural anti-holomorphic isomorphism
between $\mathscr L_1$ and $\mathscr L_2$.
This isomorphism induces an anti-holomorphic involution on the total space of the bundle $\mathscr L_1\oplus\mathscr L_2$.
Moreover, since all the planes $H_1,\cdots,H_n$ have to be real from the beginning, the variety $\mathscr Z$ is real under the above anti-holomorphic involution.

Thus we have obtained the following situation:
\begin{equation}
\label{cd1}
 \CD
 @.\mathscr Z@>{\text{inclusion}}>>\mathscr L_1\oplus\mathscr L_2\\
@.@V{q}VV @VVV \\
\mathbb{CP}^3@<{p}<<\mathscr Q\backslash\{\bm o\}@=\mathscr Q\backslash\{\bm o\}@>{\pi}>>\mathbb C.\\
 \endCD
 \end{equation}

Then we have the following

\begin{theorem}\label{thm:LBGH}
For $s\in\mathbb C$, put $\,Z_s:=(\pi\circ q)^{-1}(s)$. 
Then we have the following.
(i) If $s=1$, then $Z_s$ is biholomorphic to a Zariski open subset of the projective model \eqref{LB-n} of the LeBrun metric we have taken (just after Lemma \ref{lemma:rest4}).
(ii) If $s\in\mathbb R^{\times}$ and $|s|\le 1$, then $Z_s$ is biholomorphic to a Zariski open subset of a projective model  of a LeBrun twistor space on $n \mathbb{CP}^2$,
(iii) If $s=0$, then $Z_s$ is biholomorphic to the variety $\tilde Z$ of  \eqref{Htw}.
(iv) In (i) -(iii) the biholomorphic maps can be taken in such a way that the maps preserve the real structure.
\end{theorem}

Thus we have explicitly given a degeneration which connects any LeBrun metric on $n\CP^2$ and 
a Gibbons-Hawking metric.
(Note that for $s\in\mathbb R$ with sufficiently large $|s|$ (i.e.\,backward degeneration), $Z_s$ is {\em not} a projective model of a LeBrun twistor space.)

\begin{proof}
If $s\in\mathbb R\backslash\{0\}$, the restriction $h_i|_{Q_s}$ is a real section of $\mathscr O(1,1)$.
When $s=1$, by the choice of $h_i$, this defines the discriminant curve $\mathscr C_i$.
We write $h_i=x_0-b_ix_1-c_ix_2-\ol c_ix_3$.
Then by Lemma \ref{lemma:hi}, we have $|b_i|^2+4|c_i|^2<1/1^2=1$.
Therefore, if $s\in\mathbb R$ satisfies $|s|\le 1$, the inequality $|b_i|^2+4|c_i|^2<1/s^2$ holds.
Hence the plane section $Q_s\cap H_i$ still corresponds to a point on $\mathscr H^3$.
So the restriction  of the right-hand side $h_1\cdots h_n$ of \eqref{eqn:scrZ} to $Q_s$ defines (mutually distinct) $n$ curves of bidegree $(1,1)$ corresponding to points on $\mathscr H^3$.

On the other hand, for the left-hand side, by Lemma \ref{lemma:rest4}, the restrictions of $\mathscr L_1$ and $\mathscr L_2$ to the fiber $Q_s$ are isomorphic to $\mathscr O(n-1,1)$ and  $\mathscr O(1,n-1)$ respectively.
Therefore the defining equation of $Z_s$ is the same as the defining equation \eqref{LB-n}, where we read the equation as defined on the Zariski open subset $\mathscr O(1,n-1)\oplus\mathscr O(n-1,1)$ of the projective plane bundle $\mathbb P(\mathscr O(1,n-1)\oplus\mathscr O(n-1,1)\oplus \mathscr O)$.
Hence we obtain the claims (i) and (ii).

If $s=0$, then again by Lemma \ref{lemma:rest4}, both of the restrictions $x|_{Q_0}$ and $y|_{Q_0}$ belong to $\mathscr O([n\bm e])$, where $\bm e$ is an edge of the cone as before.
Moreover if $r:Q_0\to\mathbb{CP}^1$ denotes the projection to the conic, then $\mathscr O([\bm e])$ is isomorphic to $r^*\mathscr O(1)$.
Therefore 
$x|_{Q_0}$ and $y|_{Q_0}$ belong to $r^*\mathscr O(n)$.
Hence over $Q_0$, the point $(x,y)$ belongs to the bundle 
$r^*\mathscr O(n)^{\oplus 2}$.
The total space of this bundle is isomorphic to that  of $\mathscr O(n)^{\oplus 2}\oplus\mathscr O(2)$.
On the other hand,  under the identification $Q_0\backslash\{\bm{o}\}\simeq\mathscr O(2)$ the restrictions $h_i|_{Q_0}$ is clearly a section of $\mathscr O(2)$.
Further, by the last remark in Lemma \ref{lemma:hi}, the zeros of $h_i|_{Q_0}$ does not
hit the vertex $\bm o$.
Hence if $z$ denotes a coordinate on fibers of $\mathscr O(2)$,
then  $h_i|_{Q_0}$ can be written as $z-p_i(u)$ for some quadratic polynomial $p_i(u)$.

We show that the coefficients of $p_i(u)$ are subject to the reality conditions in Hitchin's case.
To see this, by the equation $\{x_1^2+x_2x_3=0\}$ of the conic,
the ratio $x_1/x_2$ can be taken as a non-homogeneous coordinate on the conic.
So we can put $u=x_1/x_2$.
On the other hand, the plane $H_i$ is defined by $x_0-b_ix_1-c_ix_2-\ol c_ix_3=0$ with $b_i\in\mathbb R$ and $c_i\in\mathbb C$.
Dividing the equation by $x_2$, noting that $x_0/x_2$ can be used as a fiber coordinate on $\mathscr O(2)$, and also the relation $u^2+(x_3/x_2)=0$ valid on the cone,
we obtain the equation $z-b_iu-c_i-\ol c_i(-u^2)=0$.
This means 
$$
p_i(u)=-\ol c_iu^2+b_iu+c_i.
$$
Therefore we obtained the required relations for the coefficients, and the defining equation of $Q_0$ is written exactly in the form \eqref{Htw}.
Hence we obtain (iii).

Finally, by the choice of our real structure on the bundle $\mathscr L_1\oplus\mathscr L_2$, 
on $Z_s$ with  $s\in\mathbb R^{\times}$, it agrees with that on $\mathscr O(n-1,1)\oplus\mathscr O(1,n-1)$ explained in the beginning of Section 2.
This means that the isomorphisms in (i) and (ii) preserve the real structure.
By a similar reason, the real structures also agree on the central fiber $Z_0$, meaning (iv).
\end{proof} 

From the proof, the Zariski open subset of the projective model $X$ to which $Z_s$ $(s\in\mathbb R^{\times}$ and $|s|\le 1)$ is biholomorphic, is exactly the subset
\begin{align}\label{Z-open}
\{(x,y,z)\in X\set z\neq0\}.
\end{align}
In the notation of Section 2, the complement of the set \eqref{Z-open}  in $X$ is exactly the two divisors $E$ and $\ol E$.
Recalling that in the construction of the LeBrun twistor space, $E$ and $\ol E$ are blown-down to $\CP^1$.
These curves clearly form a conjugate pair, so we denote them by $C_1$ and $\ol C_1$.
Then with respect to the twistor fibration to $n \mathbb{CP}^2$,  $C_1$ and $\ol C_1$ are exactly over the 2-sphere fixed by the semi-free U(1)-action on $n \mathbb{CP}^2$.
Thus by taking a small resolution of double points into account,  it is possible to say that {\em the complement of $C_1\cup\ol C_1$ in the LeBrun twistor spaces on $n \mathbb{CP}^2$ can be deformed to the Hitchin's twistor space of the minimal resolution of $\mathbb C^2/\Gamma$. }
If $\pi:Z_{LB}\to n \mathbb{CP}^2$ denotes the twistor fibration map, then we can alternatively say that {\em the non-Zariski open subset $\pi^{-1}(n \mathbb{CP}^2\backslash\pi(C_1))$ can be deformed to the Hitchin's twistor space.}
This means that the non-holomorphic submanifold $\pi^{-1}(\pi(C_1))$ (which are union of twistor lines parameterized by $\pi(C_1)\simeq S^2$) shrinks to the single twistor line over the orbifold point.
In particular, a non-holomorphic submanifold shrinks to a holomorphic submanifold.
On the other hand, since the locus $\pi(C_1)\,(\simeq S^2)$ is canonically identified with the boundary of $\mathscr H^3$, from the metric viewpoint, this is entirely reasonable.

For explaining what happens for the metrics in this degeneration, we first recall that if $s\in\mathbb R^{\times}$ and $|s|\le 1$, the base space  $Q=Q_s$ of the conic bundle $Z_s\to Q_s$ (the restriction of $q:\mathscr Z\to\mathscr Q$) is exactly the minitwistor space of the upper-half space $\mathscr H^3$ (regarded as an Einstein-Weyl space).
By Lemma \ref{lemma:hi},
the space $\mathscr H^3$ is identified with the ellipsoid 
$$
\mathscr B(s):=\{(b,c)\in\mathbb R\times\mathbb C\set
b^2+4|c|^2<s^{-2}\}.
$$
Then clearly we have $\lim_{s\to 0}\mathscr B(s)=\mathbb R^3$.
On the other hand, when we defined the family $q:\mathscr Z\to \mathscr Q$,  we had {\em fixed} the $n$ planes $H_1,\cdots,H_n$.
Let  $p_1,\cdots,p_n$ be the corresponding monopole points respectively. 
These $n$ points naturally belong to the ellipsoid $\mathscr B(s)$.
Further, as the right-hand side of  the equation in  \eqref{eqn:scrZ} is independent of $s$, these $n$ points in $\mathscr B(s)$ do not move.
Thus we have seen that {\em when $s$ goes to zero, while the Einstein-Weyl space $\mathscr B(s)$ tends to the Euclidean space $\mathbb R^3$, the monopole points stay fixed}.

Next we show that {\em if we pull-back the family $q:\mathscr Z\to\mathscr Q$ to $Q\times\mathbb C^*$ by the trivializing map $\phi$ (see \eqref{isom15}), then on the central fiber there appears the twistor space of the LeBrun metric on $\mathscr O(-n)$}.
For this,  if we put $u=x_2/(sx_0-x_1)$ and $v=x_3/(sx_0-x_1)$, then by a simple computation, on the quadric $Q_s$ $(s\neq 0)$, 
we have
\begin{align}\label{pQs}
(x_0,x_1,x_2,x_3)=(uv+1, s(uv-1),2su,2sv).
\end{align}
On the other hand, by Lemma \ref{lemma:hi},  the plane $H_i$ is defined by the equation \eqref{H_i} with the coefficients $b, c$ replaced by $b_i,c_i$, where $(b_i,c_i)\in\mathbb R\times\mathbb C$ satisfies the inequality \eqref{ellipsoid}.
Substituting \eqref{pQs} to the equation of $H_i$, we get
\begin{equation}\label{pH_i}
(1-b_is)uv-2sc_iu-2s\ol c_iv + (1+b_is)=0.
\end{equation}
Let $\tilde h_i=\tilde h_i(u,v,s)$ be the left-hand side of \eqref{pH_i}.
Then $\tilde h_i$ is a defining equation of the divisor $\phi^{-1}(H_i)$ in $\mathbb{CP}^1\times\mathbb{CP}^1\times\mathbb C$.
Therefore,  the pull-back family $\phi^*\mathscr Z \to
\mathbb{CP}^1\times\mathbb{CP}^1\times\mathbb C$ is given by the equation
\begin{align}
xy=\tilde h_1\tilde h_2\cdots\tilde h_n,
\end{align}
where $(x,y)\in\mathscr O(n-1,1)\oplus\mathscr O(1,n-1)$ by Lemma \ref{lemma:rest4} as before.
 Then by \eqref{pH_i}, the defining equation of the inverse image of $\mathbb{CP}^1\times\CP^1\times\{0\}$ (under the projection $ \phi^*\mathscr Z \to
\mathbb{CP}^1\times\mathbb{CP}^1\times\mathbb C$) is given by
\begin{align}
xy= (uv+1)^n.
\end{align}
By Theorem \ref{thm:O(-n)}, this is exactly the defining equation of (a Zariski open subset of)
the twistor space of the LeBrun metric on  $\mathscr O(-n)$, as required.

If one wants not a subset but the entire twistor space, 
it is enough to compactify the fibers of the family $\mathscr Z(\to\mathscr Q)\to\mathbb C$ by 
considering the projective bundle $\mathbb P(\mathscr L_1\oplus\mathscr L_2\oplus\mathscr O)$ instead of $\mathscr L_1\oplus\mathscr L_2$, and pull it back by $\phi$.

\section{Discussions about degenerations of the metrics}

So far we have studied degenerations of twistor spaces.
But originally  these are of course motivated by understanding degenerations of (anti)-self-dual metrics or conformal classes on the 4-manifolds.
In this final section we  briefly discuss what happens for base 4-manifolds, by picking up 
 K\"ahler representatives of the conformal classes on (dense) open subsets of the 4-manifolds.
For a  thorough investigation in this direction, we refer a paper by Viaclovsky \cite{V}.

We recall that LeBrun \cite{LB91} first constructs a scalar-flat K\"ahler metric on a collinear blow-up of $\mathbb C^2$ explicitly, and then verifies that the metric can be extended to a one point compactification after adjusting a conformal gauge.
Here for avoiding confusion we call these scalar flat K\"ahler metrics as LeBrun {\em metrics}.

As showed by Tian-Viaclovsky \cite[Theorem 1.1]{TV05b},
any sequence of complete K\"ahler   metrics with constant scalar curvature on a 4-manifold satisfying some boundedness condition on the curvature tensor and Sobolev constants has a converging subsequence with the limit being a K\"ahler orbifold, where the convergence is in the Gromov-Hausdorff  sense as {\em pointed}   spaces. 
For avoiding confusion with monopole centers, we refer the point of pointed spaces as a {\em base point}.

Let $\{\{p_{i1},p_{i2},\cdots,p_{in}\}\subset \mathscr H^3\set i=1,2,\cdots\}$ be a sequence of $n$ points on the hyperbolic space, and $\{g_i\}$ the sequence of the LeBrun metrics which have these $n$ points as the set of monopole points.
First as in Section 2  we 
consider the situation that  the first point $p_{i1}$ goes to infinity as $i\to \infty$ while 
the remaining $(n-1)$ points stay fixed.
Then if we take a compact domain $K\subset\mathscr H^3$ and if we choose as a base point a point $q_i$ which belongs to the subset over $K$, then since the Green function centered at $p_{i1}$ tends to zero as $i\to \infty$ on the $R$-ball centered at $q_i$ for any $R>0$,
the space with LeBrun metric $g_i$ converges to the space equipped with the LeBrun metric having the fixed $(n-1)$ points as the set of monopole points.
The degeneration from $n\mathbb{CP}^2$ to $(n-1)\mathbb{CP}^2$ constructed in Section 2 is 
the twistor translation of this degeneration.

Next consider the same sequence of $n$ points as in the previous paragraph but this time as base points we choose the moving point $p_{i1}$, viewed as a point on $n\mathbb{CP}^2$.
Then since the distance between the fixed $(n-1)$ points and the base point goes to infinity as $i\to \infty$, the limit as pointed space should be the space equipped with the LeBrun metric with one monopole point; namely the Burns metric.
This is rather the space which disappeared in the last degeneration, and
thus changing a base point can give a different limit.

Also  it is natural to expect that in the degeneration of the twistor spaces constructed in Section 2, 
at the limit (central fiber), we  find not only a twistor space of $(n-1)\mathbb{CP}^2$ but also the twistor space of Fubini-Study metric (which is a conformal compactification of the Burns metric).
However, even if we regard our degeneration as a family $\mathscr Z\to \mathbb C(\lambda)$ and try to perform the birational transformations within the 4-fold $\mathscr Z$, I could not find the flag twistor space at the central fiber.
We note that the existence of such a model is guaranteed by the framework of Donaldson-Friedman \cite{DF89}.
If one can find the flag twistor space at the central fiber of the family, it provides an explicit realization of  the Donaldson-Friedman model for the case of the present degeneration,
not relying on deformation theory of complex spaces.

Next in accordance with the degeneration taken up in Section 3 we consider the sequence $\{\{p_{i1},p_{i2},\cdots,p_{in}\}\subset \mathscr H^3\set i=1,2,\cdots\}$ of $n$ points for which all $p_{ij}$-s approach to a point $p\in\mathscr H^3$ as $i\to \infty$.
In this case, as was explicitly shown in coordinates by LeBrun \cite[page 235-236]{LB91}, the limit is a space equipped with the LeBrun's scalar flat K\"ahler metric on $\mathscr O(-n)$.
However we note that we have to choose a conformal gauge which is different from the above LeBrun metrics on the collinear blowup of $\mathbb C^2$; instead another K\"ahler representative found in \cite[page 243--244]{LB91} has to be chosen to get the above limit.

As was discovered in \cite[page 236--237]{LB91}, if we look at this degeneration more carefully, it provides a typical example of bubbling off phenomena of ALE spaces.
Namely, when all the monopole points $p_{i1},\cdots, p_{in}$ of LeBrun metrics become closer to a point  $p\in \mathscr H^3$ in such a way that their angular positions relative to $p$ as well as the ratio of the distances from $p$ are preserved, around the point $p$ we make a sequence of rescalings in a way that the distances from $p$ become constant.
Then at the limit the curvature of the hyperbolic space becomes zero and as a result a Gibbons-Hawking space with $n$ monopole points on $\mathbb R^3$ bubble off from the point $p$.

To understand this bubble off through our degeneration of the twistor spaces constructed in Section 4, 
for each $s\in\mathbb R$ with $0<s<1$, we define a dilation $\psi_s:\mathbb R^3\to \mathbb R^3$
by $\psi_s(b,c)=(b/s,c/s)$ for $(b,c)\in\mathbb R\times\mathbb C=\mathbb R^3$. Then for the ellipsoids  we have $\psi_s(\mathscr B(1))=\mathscr B(s)$.
Hence pulling back by $\psi_s$, $\mathscr B(s)$ becomes $\mathscr B(1)$ and the point 
$p_j=(b_j,c_j)\in\mathscr B(s)$ is pulled back to $s. p_j=(sb_j,sc_j)$.
Thus by letting $s\to 0$,  all the monopole points become closer to the origin,
while their relative position is exactly as explained in the last paragraph.

Finally it seems natural to expect that, when any sequence $\{\{p_{i1},p_{i2},\cdots,p_{in}\}\subset \mathscr H^3\set i=1,2,\cdots\}$ of $n$ points is given, possible limits (as pointed spaces) of the associated LeBrun metrics, implied by the above Tian-Viaclovsky's convergence theorem, are  LeBrun orbifolds with $k$ ($0\le k \le n$) monopole points, where the LeBrun orbifold metric means a natural generalization of the LeBrun metrics, allowing some of the monopole points to coincide (see \cite[Section 2.2]{V} for the precise definition).
Then it would be possible to say that the degenerations considered in Sections 2 and 3 are two extremal cases among these degenerations.
\vspace{3mm}

\end{document}